\documentclass[12pt]{amsart}

\usepackage{amsmath}

\usepackage{amscd}

\usepackage{graphicx}

\usepackage{latexsym}

\usepackage{hyperref}

\usepackage{times}

\oddsidemargin .25in \theoremstyle{plain}

\newtheorem{Thm}{Theorem}

\newtheorem{Cor}[Thm]{Corollary}

\newtheorem{Prop}[Thm]{Proposition}

\newtheorem{Def}[Thm]{Definition}

\newtheorem{Rem}[Thm]{Remark}

\newcommand{\OB}{\mathfrak{{ob}}}

\begin{document}

\title[]{Explicit horizontal open books on some Seifert fibered 3--manifolds}

\author{Burak Ozbagci}

\address{Department of Mathematics \\ Ko\c{c} University \\ Istanbul, Turkey}

\email{bozbagci@ku.edu.tr}

\address{ Current: School of Mathematics, Georgia Institute
of Technology, Atlanta,
Georgia, USA}

\email{bozbagci@math.gatech.edu}


\date{\today}

\thanks{The author was partially supported by the Turkish Academy of Sciences and by
the NSF Focused Research Grant FRG-024466.}

\begin{abstract}

We describe explicit horizontal open books on some Seifert fibered
3--manifolds. We show that the contact structures compatible with
these horizontal open books are Stein fillable and horizontal as well. 
Moreover we draw surgery diagrams for some of these contact structures.
\end{abstract}

\maketitle

\section{Introduction}

An open book on a Seifert fibered 3--manifold is called
\emph{horizontal} if its binding is a collection of some fibers and
its pages are positively transverse to the Seifert fibration. Here
we require that the orientation induced on the binding by the pages
coincides with the orientation of the fibers induced by the
fibration. In this article, we construct explicit horizontal open
books on some Seifert fibered 3--manifolds.  We show that the
monodromies of these horizontal open books are given by products of
(multiple) right-handed Dehn twists along boundary parallel curves.
Consequently, by a theorem of Giroux \cite{gi}, the contact
structures compatible with our horizontal open books are Stein
fillable. We also show that these contact structures are horizontal,
i.e., the contact planes are positively transverse to the fibers of
the Seifert fibrations. Finally we draw surgery diagrams for those
contact structures which are compatible with planar open books.

It is well-known that any Seifert fibered 3--manifold can be
obtained by a plumbing of circle bundles according to a ``standard''
diagram. The standard diagram is a star  shaped graph, with $k$
linear branches, whose central vertex is a circle bundle over an
arbitrary closed surface with arbitrary Euler number $n$ and all the
other vertices are circle bundles over $S^2$ with Euler numbers less
than or equal to $-2$. In this language a vertex always represents a
circle bundle and an edge between two vertices represents plumbing
of the circle bundles represented by the vertices. We will call a
standard diagram with $n+k \leq 0$ as \emph{non-positive} standard
plumbing diagram.

In this article we consider Seifert fibered 3--manifolds given by
plumbing of circle bundles according to a star shaped graph, with a
central vertex, a circle bundle over an arbitrary closed surface
with non-positive Euler number, and $k$ other vertices which are
circle bundles over $S^2$ with positive Euler numbers, all connected
to the central vertex. Note that these are a priori non-standard
plumbing diagrams. We construct open books on these Seifert fibered
3--manifolds which are horizontal with respect to their Seifert
fibrations. However, we can convert these diagrams into non-positive
standard plumbing diagrams so that the Euler numbers of all the
vertices except the central vertex are equal to $-2$. In \cite{eo},
explicit open books were constructed on Seifert fibered 3--manifolds
which can be described by non-positive standard plumbing diagrams.
Moreover it was shown that these open books are horizontal with
respect to the circle bundles involved in the plumbing descriptions
of the Seifert fibered 3--manifolds at hand. It turns out that the
horizontal open book we construct in this article using a 
non-standard plumbing description of a Seifert fibered 3-manifold is
isomorphic to the horizontal open book constructed  in \cite{eo} 
starting with the non-positive 
standard plumbing description of the same manifold.
The difference is that here horizontal means transverse to the
Seifert fibers and in \cite{eo} horizontal means transverse to the
circle bundles involved in a plumbing description of a Seifert
fibered 3--manifold. It is not clear to the author whether or not
these two different definitions of being horizontal are equivalent
in general.

In particular suppose that the circle bundle at the central vertex
in the plumbings are also over $S^2$. Then the open books we
construct are planar, which means that a page has genus zero. The
significance of finding contact structures compatible with planar
open books stems from the recent work of Abbas, Cieliebak and Hofer
\cite{ach} who proved the Weinstein conjecture for those contact
structures. In \cite{sc}, Sch\"{o}nenberger also constructs planar
open books compatible with Stein fillable contact structures which
are given by Legendrian surgeries on any Legendrian realization of a
non-positive standard diagram of a Seifert fibered 3--manifold. We
can pin down the Stein fillable horizontal contact structures
compatible with the horizontal open books we construct, by comparing
our construction with Sch\"{o}nenberger's method.

\vspace{1ex}

\noindent {\bf {Acknowledgement}}: We are grateful to John Etnyre
for helpful conversations and Andr\'{a}s Stipsicz  for generously
sharing his expertise and commenting on the first draft of this paper.

\section{Open book decompositions and contact structures}

We will assume throughout this paper that a contact structure
$\xi=\ker \alpha$ is coorientable (i.e., $\alpha$ is a global
1--form) and positive (i.e., $\alpha \wedge d\alpha >0 $ ). In the
following we describe the compatibility of an open book
decomposition with a given contact structure on a 3--manifold.

Suppose that for an oriented link $L$ in a closed and oriented
3--manifold $Y$ the complement $Y\setminus L$ fibers over the
circle as $\pi \colon Y \setminus L \to S^1$ such that
$\pi^{-1}(\theta) = \Sigma_\theta $ is the interior of a compact
surface bounding $L$, for all $\theta \in S^1$. Then $(L, \pi)$ is
called an \emph{open book decomposition} (or just an \emph{open
book}) of $Y$. For each $\theta \in S^1$, the surface
$\Sigma_\theta$ is called a \emph{page}, while $L$ the
\emph{binding} of the open book. The monodromy of the fibration
$\pi$ is defined as the diffeomorphism of a fixed page which is
given by the first return map of a flow that is transverse to the
pages and meridional near the binding. The isotopy class of this
diffeomorphism is independent of the chosen flow and we will refer
to that as the \emph{monodromy} of the open book decomposition.

\vspace{1ex}

\begin{figure}[ht]

  \begin{center}

     \includegraphics{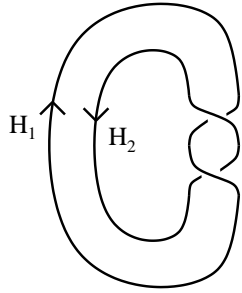}

   \caption{Positive Hopf link $H$ in $S^3$ }

  \label{positive-hopf}

    \end{center}

  \end{figure}

\begin{figure}[ht]

  \begin{center}

     \includegraphics{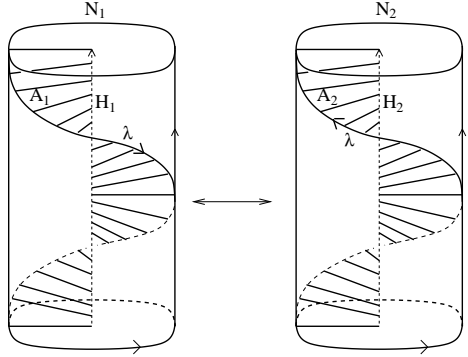}

   \caption{By gluing the two solid tori we can visualize an open book of
   $S^3$, with binding $H$, whose page is obtained by gluing $A_1$ with $A_2$ along $\lambda$.}

  \label{heegaard}

    \end{center}

  \end{figure}

\noindent \textbf{Example}: Consider the positive Hopf link $H$
(see Figure~\ref{positive-hopf}) in $S^3$. Note that $S^3$ has a
Heegaard splitting of genus one, and each piece in the splitting
can be taken as a regular solid torus neighborhood $N_i$ (cf.
Figure~\ref{heegaard}) of the component $ H_i$ of  $H$, for $ i =
1,2$.

Take the core circle $H_1$ of $N_1$ and consider the annulus $A_1$
bounded by $H_1$ and a $(1, -1)$-curve $\lambda$ on $\partial
N_1$. Observe that the $\lambda$ becomes a $(-1, 1)$-curve on
$\partial N_2$. There is another annulus $A_2$ bounded by $H_2$
and $-\lambda$ on $\partial N_2$. Then the annulus $A = A_1 \cup
A_2$ is bounded by  the Hopf link $H$ in $S^3 = N_1 \cup N_2$. Now by
foliating each $N_i$ by (a circle worth of) annuli $A_i$ and
gluing the corresponding leafs (cf. Figure~\ref{heegaard}) we get
a fibration $S^3 \setminus H \rightarrow S^1$ of the complement of
$H$ in $ S^3$  with annuli pages $A$ each of which is bounded by $H$.
Note that the annulus $A$ is a copy of the obvious Seifert surface
bounded by  $H= H_1 \cup H_2$  in Figure~\ref{positive-hopf}.

To see the monodromy of the resulting open book take a vector
field in $N_i$ so that it has constant slope $1$ on $\partial
N_i$, and its slope smoothly decreases to zero as we move towards
$H_i$ on any horizontal ray so that it becomes meridional near the
core circle $H_i$. Note that this defines a vector field on $S^3=
N_1 \cup N_2$ which is transverse to the pages of the open book
and becomes meridional near the binding $H$. The first return map
of the flow induced by this vector field is a right-handed Dehn
twist along the center circle of a fixed page $A$ in the open
book.

\vspace{1ex}

Every closed and oriented 3--manifold admits an open book
decomposition (cf. \cite{al}) as well as a contact structure (cf.
\cite{marti}). An open book decomposition and a contact structure
is related by the following definition.

{\Def \label{compatible} An open book decomposition $(L,\pi)$ of a
3--manifold $Y$ and a contact structure $\xi$ on $Y$ are called
\emph{compatible} if $\xi$ can be represented by a contact form
$\alpha$ such that $\alpha ( L) > 0$ and $d \alpha > 0$ on every
page.}

\vspace{1ex}

\noindent \textbf{Example}: Consider the unit sphere $S^3 = \{
(z_1,z_2) : {\vert z_1 \vert }^2 + {\vert z_2 \vert }^2 =1 \}
\subset \mathbb{C}^2$. In these coordinates the positive Hopf link
is given by $$H= \{(z_1,z_2) :z_1 z_2 = 0 \}= H_1 \cup  H_2 $$
where $H_1 = \{(z_1,z_2) : z_2 = 0 \}$ and $H_2= \{(z_1,z_2) : z_1
= 0 \}$. Let $\pi : S^3 \setminus H \to S^1$ be defined as
$$\pi(z_1, z_2)= \frac{z_1 z_2}{\vert z_1 z_2 \vert} \in S^1
\subset \mathbb{C}.$$ In polar coordinates the map $\pi$ can be
expressed as : $\pi (r_1 e^{i\theta_1} ,r_2 e^{i\theta_2} ) =
e^{i(\theta_1 + \theta_2)}.$ Then for a fixed $\theta_0 \in S^1$
we have $$\pi^{-1} (\theta_0) = \{( \theta_1, r_1,  \theta_2, r_2)
: \theta_1 + \theta_2=\theta_0 , r_1^2 + r_2^2 = 1 \}.$$ Note that
the equation $r_1^2 + r_2^2 = 1$ gives an arc between a point in
$H_1$ with coordinate $\theta_1$ and the point in $H_2$ with
coordinate $\theta_2 = \theta_0 - \theta_1$. The union of these
arcs will sweep out the annulus $\pi^{-1} (\theta_0)$ which bounds
$H$. On the other hand, the standard contact structure $\xi$ in
$S^3$ can be given as the kernel of the 1--form $\alpha = r_1^2 d
\theta_1 + r_2^2 d \theta_2,$ restricted to $S^3 \subset
\mathbb{C}^2$. To show that $\alpha$ is compatible with the open
book above we first observe that $\alpha (H_i)= d\theta_i
(\partial_{\theta_i}) = 1$. Next we parameterize $\pi^{-1}
(\theta_0)$ by $$f(r,\theta)= (r, \theta, \sqrt{1-r^2}, \theta_0 -
\theta_1).$$ Thus $f^*(\alpha)= (2r^2 -1)d\theta$ and
$f^*(d\alpha)= df^*(\alpha)= 4rdrd\theta$. We conclude that
$d\alpha$ is an area form on $\pi^{-1} (\theta_0)$, for all
$\theta_0 \in S^1$.

\vspace{1ex}

We have the following fundamental results at our disposal to study the
contact topology of 3--manifolds.

{\Thm [Thurston-Winkelnkemper \cite{tw}] \label{tw} Every open book admits a
compatible contact structure.}

\vspace{1ex}

Conversely,

{\Thm [Giroux \cite{gi}] \label{giroux} Every contact 3--manifold
admits a compatible open book.
Moreover two contact structures
compatible with the same open book are isotopic.}

\vspace{1ex}

We refer the reader to \cite{et} and \cite{ozst} for more on the
correspondence between open books and contact structures.

\section{Explicit construction of horizontal
open books on some Seifert fibered 3--manifolds} \label{construction}

\begin{figure}[ht]

  \begin{center}

     \includegraphics{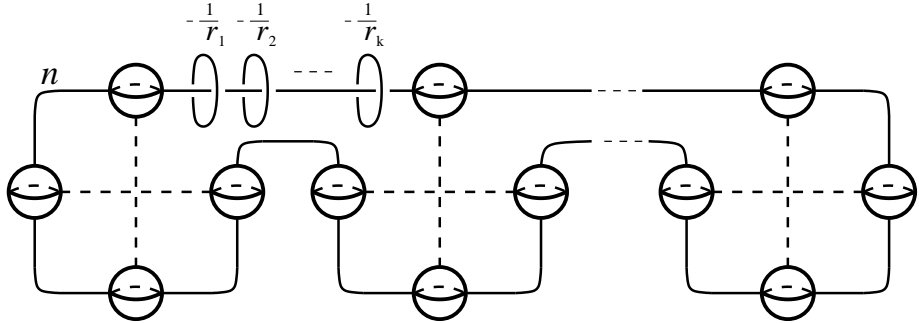}

   \caption{Seifert fibered
3--manifold with Seifert invariants $(g,n; r_1, \ldots , r_k)$.
(There are $2g$ 1--handles in the figure.)}

  \label{seifert-fibered}

    \end{center}

  \end{figure}

Consider the $S^1$--bundle $Y_{g,n}$ over a closed and oriented
genus $g$ surface $\Sigma$ with Euler number $e(Y_{n,g})=n $, as the
boundary of the corresponding $D^2$--bundle over $\Sigma$ with one 0--handle,
$2g$ 1--handles and one 2--handle (with framing $n$) as in
Figure~\ref{seifert-fibered}.  The 3--manifold $Y$ obtained from
$Y_{g,n}$ by performing $-1/r_1, \ldots, -1/r_k$ surgeries on $k$
distinct fibers (cf. Figure~\ref{seifert-fibered}) is called a
\emph{Seifert fibered 3--manifold} with Seifert invariants $(g,n;
r_1, \ldots , r_k)$ where $r_i \in \mathbb{Q} $. Notice that
according to this convention the surgery coefficients are negative
reciprocals of the given data.

An open book on a Seifert fibered 3--manifold is called
\emph{horizontal} if its binding is a collection of some fibers
and the interiors of its pages are positively transverse to the
Seifert fibration. Here we require that the orientation induced on
the binding by the pages coincides with the orientation of the
fibers induced by the fibration. In this article we construct
horizontal open books on Seifert fibered 3--manifolds with
invariants $(g,n; r_1, \ldots , r_k)$, where $g \geq 0$ is
arbitrary and
$$n \leq 0 < -\frac{1}{ r_i} \in \mathbb{N}, \; \mbox{for}\;  1
\leq i \leq k .$$

Let $\Sigma$ denote a closed oriented surface and let $M$ denote the
trivial circle bundle over $\Sigma$, i.e., $M = S^1 \times \Sigma$.
Perform a $+1$--surgery along circle fiber $K$ of $M$. The resulting
3--manifold $M^\prime$ is a circle bundle over $\Sigma$ whose Euler
number is equal to $-1$. In \cite{eo} we described an explicit
horizontal open book on $M^\prime$: A fiber becomes the binding,
every page is a once punctured $\Sigma$, and the monodromy is
right-handed Dehn twist along a curve parallel to the binding. Then
as discussed in \cite{eo}, to obtain a horizontal open book on a
circle bundle with Euler number $ n < 0 $ we need to perform
$+1$--surgeries along $\vert n \vert$ distinct fibers.

We claim that we can generalize the arguments used in \cite{eo},
for a $+1$--surgery,  to a $p$--surgery (for any  $0 < p \in
\mathbb{N}$) to obtain horizontal open books on some Seifert
fibered 3--manifolds. To perform a $p$--surgery on a fiber $K$ of
$M$ we first remove a solid torus neighborhood $N= K \times D^2$
of $K$ from $M$. Note that by removing $N$ from $M$ we puncture
once each $\Sigma$ in $M = S^1 \times \Sigma$ to get $S^1 \times
\widetilde{\Sigma}$, where $\widetilde{\Sigma}= \Sigma \setminus
D^2$. Now we will glue a solid torus back to $S^1 \times
\widetilde{\Sigma}$ along its boundary torus $S^1 \times
\partial \widetilde{\Sigma}$ in order to perform our surgery.

\begin{figure}[ht]

  \begin{center}

     \includegraphics{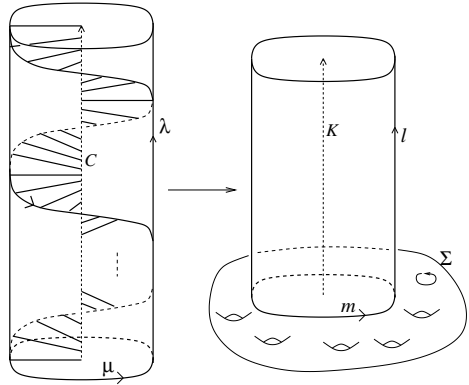}

   \caption{The surgery torus is obtained by
 identifying the top and the bottom of the cylinder shown on
the left-hand side. The annulus depicted inside the torus is
wrapped $p$-times around the core circle in a right-handed manner.
On the right hand side we depict the complement $S^1 \times
\widetilde{\Sigma}$ obtained by removing a neighborhood of a
circle fiber from $S^1 \times \Sigma$. }
\label{surgerytorus-multiple}

    \end{center}

  \end{figure}

In the discussion below we will make no distinction
between curves on a surface and the homology classes they represent, to
simplify the notation.
Consider the solid torus $S^1 \times D^2$ shown on the left-hand
side in Figure~\ref{surgerytorus-multiple}. Let $\mu$ and
$\lambda$ be the meridian and the longitude pair of $S^1 \times
\partial D^2$. Let $m$ and $l$ denote the meridian and the
longitude pair in the boundary of $S^1 \times \widetilde{\Sigma}$.
Note that the base surface is oriented (as depicted in
Figure~\ref{surgerytorus-multiple}) and the orientation induced on
$\partial  \widetilde{\Sigma}$ is the opposite of the orientation
of $m$.  We glue the solid torus $S^1 \times D^2$ to $S^1 \times
\widetilde{\Sigma}$ by an \textit{orientation preserving}
diffeomorphism $ S^1 \times
\partial D^2 \to S^1 \times \partial \widetilde{\Sigma}$
which sends $\mu$ to $pm+l$ and $\lambda$ to $ (p^2-1)m+pl$. (This
map is in fact orientation reversing if we use the orientation on
the torus  $ S^1 \times \partial \widetilde{\Sigma}$ which is
induced from $ S^1 \times \widetilde{\Sigma}$. ) The resulting
3--manifold $M^\prime$ will be oriented extending the orientation
on $S^1 \times \widetilde{\Sigma}$ induced from $M$.

We also depict, in Figure~\ref{surgerytorus-multiple},  a leaf (an
annulus) of a foliation on the solid torus $S^1 \times D^2$ that
we will glue to perform $p$--surgery. Note that under the surgery
map $ p \lambda - (p^2 -1) \mu $ is mapped onto $l$ which
parametrizes the $S^1$ factor in $S^1 \times \widetilde{\Sigma}$.
We can view $ S^1 \times (D^2 - \{ 0\})$ as a disjoint union of
concentric tori and foliate each torus by (a circle worth of)
curves isotopic to $p \lambda - (p^2 -1) \mu $. This gives a
foliation of $ S^1 \times D^2 $ into circles where the center
circle $S^1 \times \{ 0\}$ is covered $p$--times by all the other
circles in the foliation. Note that this circle foliation defines the Seifert fibration.
Observe that the Seifert fibration on
$S^1 \times D^2$ is transverse to the annuli foliation simply
because the slope of the circles is less than the slope of the
annular pages in $S^1 \times D^2$.

The boundary of a leaf consists of the core circle $C$ of $S^1
\times D^2$ and a $(p,-1)$--curve on $S^1 \times
\partial D^2$, i.e., a curve isotopic to $ p\mu - \lambda$. Each
leaf is oriented so that the induced orientation on the boundary
of a leaf is given as indicated in
Figure~\ref{surgerytorus-multiple}. The gluing diffeomorphism maps
$p\mu - \lambda$ onto $m$ so that by performing the $p$--surgery
we also glue each annulus in the foliation to a
$\widetilde{\Sigma}$ in $S^1 \times \widetilde{\Sigma}$
identifying the outer boundary component (i.e., the $(p\mu -
\lambda)$--curve) of the annulus with $\partial
\widetilde{\Sigma}$. Hence this construction yields a horizontal
open book on $M^\prime$ whose binding is $C$ (the core circle of
the surgery torus) and pages are obtained by gluing an annulus to
each $\widetilde{\Sigma}$ along $\partial \widetilde{\Sigma}$.
Notice that the pages will be oriented extending the orientation
on $\widetilde{\Sigma}$ induced from $\Sigma$. Finally we want to
point out that the core circle $C$ becomes an oriented fiber of
the Seifert fibration of $M^\prime$ over $\Sigma$. In fact $C$
becomes a singular fiber if $ p> 1$.

\begin{figure}[ht]

  \begin{center}

     \includegraphics{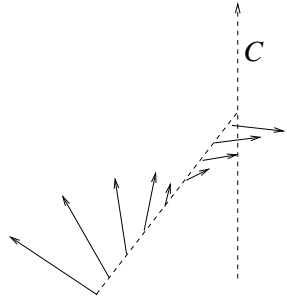}

   \caption{The vector field depicted along a ray towards
the core circle $C$ inside the surgery torus extending the
vertical one on $S^1 \times \widetilde{\Sigma}$. The vector field
becomes horizontal near the core (binding) along any ray.}
\label{flow-multiple}

    \end{center}

  \end{figure}

Next we will describe the monodromy of this open book. In order to
measure the monodromy of an open book we should choose a flow
which is transverse to the pages and meridional near the binding.
We will take a \emph{vertical} vector field pointing along the
fiber direction in $S^1 \times \widetilde{\Sigma}$ and extend it
inside the surgery torus as follows: The vertical vector field is
given by $\partial_l$ on $ S^1 \times \partial \widetilde{\Sigma}$
which is identified with $p \partial_\lambda -
(p^2-1)\partial_\mu$ on the boundary of the surgery torus. Extend
this vector field with slope $-p/(p^2 -1 )$ inside he surgery
torus (along every ray towards the core circle) by rotating
clockwise so that it becomes horizontal near the core circle as
illustrated in Figure~\ref{flow-multiple}. First observe that this
vector field is positively transverse to the pages of the open
book by construction since $- p / (p^2 -1) < -1 / p$. In other
words one can see that this vector field always points towards the
positive side of the annulus in Figure~\ref{surgerytorus-multiple}
by comparing the slope of the annulus (a page)  with the slope of
the vector field.  Next note that the first return map of the flow
will fix the points near the binding on any leaf since the vector
field is horizontal near the binding. The first return map will
fix the points on $\widetilde{\Sigma}$ as well, since the flow is
vertical (i.e., in the direction tangent to the $S^1$ factor) on
$S^1 \times \widetilde{\Sigma}$. Now take a horizontal arc (on a
leaf) connecting the core circle to the other boundary of that
leaf. Then one can see that the flow will move the points of this
arc further to  the right if we move towards the boundary. The
first return map is given by a right-handed Dehn twist along the
core circle of the leaf when we go around $\lambda$ once. But
since a circle fiber $l$ covers the core circle $\lambda$ (a
singular fiber for $p > 1$) $p$-times, the first return map of the
open book is given by the $p$-th power of a right-handed Dehn
twist along the core circle of the leaf in
Figure~\ref{surgerytorus-multiple}, which is indeed a curve
parallel to the binding of the open book.

In summary we constructed a horizontal open book on a Seifert
fibered 3--manifold with invariants $(g, 0; r)$ for $0< p=-
\frac{1}{r} \in\mathbb{N}$ by performing a $p$--surgery on a fiber
of a trivial circle bundle over a genus $g$--surface $\Sigma$.
Note that the binding is a singular fiber of the Seifert fibration
if $- \frac{1}{r}= p
> 1$.  Since surgery modifies the open book only in a neighborhood
of the surgery curve in the above construction, we can construct a
horizontal open book on the Seifert fibered 3--manifold with
invariants $(g,0; r_1, \ldots , r_k)$ where $0 < -\frac{1}{ r_i}
\in \mathbb{N} $, for $1 \leq i \leq k$. Moreover if we start with
the trivial circle bundle, perform $+1$--surgeries on $\vert n \vert$ distinct
fibers (cf. \cite{eo}) for $n < 0$, and then further perform
$-1/r_1, \ldots, -1/r_k$ surgeries on $k$ disjoint fibers,
respectively,  where $0 < -\frac{1}{ r_i} \in \mathbb{N} $, for $1
\leq i \leq k$, we will get a horizontal open book on the Seifert
fibered 3--manifold with invariants $(g,n; r_1, \ldots , r_k)$.
The binding of the resulting open book will be the union of $k+\vert n \vert$
circle fibers and a page is a genus $g$ surface with $k+\vert n \vert$
boundary components. The monodromy will be a product of (multiple)
right-handed Dehn twists along boundary parallel curves. More
precisely, there are $\vert n \vert$ boundary components each of which has one
right-handed Dehn twist and for the other $k$ components of the
boundary the exponent of the right-handed Dehn twist is given by
$p_i = -1/r_i$, for $1\leq i \leq k$. Since the monodromy is a
product of right-handed Dehn twists only, the contact structure
compatible with this open book is Stein fillable by a theorem of
Giroux (cf. \cite{gi}). Thus we showed

{\Prop \label{hor} There exists an explicit \emph{horizontal} open
book on the Seifert fibered 3--manifold with invariants $(g,n;
r_1, \ldots , r_k)$, where $n \leq  0 < -\frac{1}{ r_i} \in
\mathbb{N} $, for $1 \leq i \leq k$. The contact structure
compatible with this open book is Stein fillable.}

\vspace{1ex}

Next we consider the contact structure compatible with our open book.

{\Prop \label{transverse} The contact structure compatible with
the open book we constructed on the Seifert fibered 3--manifold
with invariants $(g,n; r_1, \ldots , r_k)$,  such that $n \leq  0
< -\frac{1}{ r_i} \in \mathbb{N} $, for $1 \leq i \leq k$, is
\emph{horizontal}, i.e., the contact planes are positively
transverse to the fibers of the Seifert fibration.}

\begin{proof}
Consider the contact structure compatible with the constructed
horizontal open book. By an isotopy of this contact structure we
may assume that the contact planes are arbitrarily close to the
tangents of the pages away from the binding (cf. Lemma 3.5 in
\cite{et}).  Since the pages of our open book are already
positively transverse to the fibers we can conclude that the
contact planes are positively transverse to the fibers away from
the binding. Note that the contact structure is still compatible
with the given open book after this isotopy.

Recall that we explicitly constructed the pages of the open book
near a component $C$ of the binding. The fibers of the circle
fibration can be viewed as straight vertical lines in the solid
cylinder on the left-hand side in
Figure~\ref{surgerytorus-multiple} before we identify the top and
the bottom. On the other hand, by the compatibility of the open
book and the contact structure, there are coordinates $(z, (r,
\theta)) $ near every component $C$, where $C$ is $\{r=0\}$ and a
page is given by setting $\theta$ equal to a constant such that
the contact structure is given by the kernel of the form $dz + r^2
d \theta$.
\begin{figure}[ht]

    \begin{center}

     \includegraphics{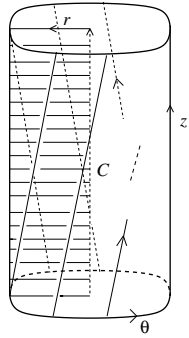}

   \caption{A page is
straightened out near a binding component $C$ in these coordinates.
A Seifert fiber is given by a $(1,p)$--curve in a neighborhood of
$C$. } \label{twist-multiple}

    \end{center}
  \end{figure}
In these coordinates the neighborhood of $C$ in
Figure~\ref{surgerytorus-multiple} is seen as in
Figure~\ref{twist-multiple}, where the annulus is straightened out
and the circle fibers are wrapped around the cylinder in a
right-handed manner. More precisely, the tangent vector to an
(oriented) fiber is given by $\partial_{\theta} + p\partial_z$.
Since
$$ (dz + r^2 d\theta) (\partial_{\theta} + p\partial_z) = r^2 +p > 0 ,$$
the contact planes are positively transverse to the circle fibers
in a neighborhood of $C$. This finishes the proof of the
proposition.
\end{proof}

{\Rem \label{hori} The proof above shows that in fact, the contact
structure compatible with a horizontal open book on any Seifert
fibered 3--manifold is horizontal. 
This is because, in a local model
near the binding, a Seifert fiber has to be of the form
$a\partial_{\theta} + b\partial_z$, for some positive integers $a$
and $b$, and indeed we have $ (dz + r^2 d\theta)
(a\partial_{\theta} + b\partial_z) = ar^2 +b> 0.$}

\vspace{1ex}

Note that the horizontal open book we constructed on the Seifert
fibered 3--manifold with invariants $(0,n; r_1, \ldots , r_k),$
where $n \leq 0 < -\frac{1}{ r_i} \in \mathbb{N}$, for $ 1\leq i
\leq k$, is planar, i.e., a page has genus zero. Recently it was
shown in \cite{ach} that the Weinstein conjecture holds for a
contact structure compatible with a planar open book.

\vspace{1ex}

The same construction of an open book will work when we perform $
p$--surgeries along fibers of a circle bundle for $p < 0$, but in
that case the orientation on the binding induced by the pages will
be the opposite of the fiber orientation. Thus the open book we
construct using $p$--surgeries will not be horizontal when $p < 0$.
Moreover we will get left-handed Dehn twists along boundary
components instead of right-handed Dehn twists in that case.
Finally, when we perform $0$--surgery along a fiber of a trivial
circle bundle over $\Sigma$, we still get an open book where the
binding is a fiber, page is once punctured $\Sigma$ and monodromy is
the identity map. Note, however, that a circle bundle does not
extend to a Seifert fibration when we do $0$--surgery on a regular
fiber. Combining all the discussion above we get the following:

{\Cor \label{boundary} Let $Y$ be a 3--manifold with an open book
whose page is a genus $g$ surface $\Sigma$ with $r$ boundary
components and let $\delta_i$ denote a curve parallel to the $i$-th
boundary component of $\Sigma$. If the monodromy of this open book
can be factorized as $\phi = \prod_{i=1}^r t_{\delta_i}^{m_i}$ for
some nonzero integers $ m_1, m_2, \cdots, m_r$ then $Y$ is a Seifert
fibered 3--manifold. This open book can be realized as a horizontal
open book with respect to the Seifert fibration if and only if
$m_i>0$ for all $i=1,2, \cdots, r$.}

\begin{proof}
The first statement and the sufficiency of the second statement
should be clear from the arguments above. To prove the necessity of
the second statement let us assume that $m_j <0 $ for some $j$. Then
one can see that there is a properly embedded arc on $\Sigma$ with
an initial point on the $j$-th component of $\partial \Sigma$ which
is not right-veering, unless $g=0$ and $r=1,2$. Hence by \cite{hkm},
this open book can not be compatible with a tight contact structure.
But it is known (cf. \cite{lm}) that any horizontal contact
structure on a Seifert fibered 3--manifold is fillable, hence tight.
This proves, by Remark~\ref{hori}, that the open book can not be
made horizontal under the given assumption. The cases we left out
above can be easily treated separately.
\end{proof}

It was shown in \cite{dg} that any contact $3$--manifold can be
obtained by contact $(\pm 1)$--surgeries along a Legendrian link in
the standard contact $S^3$. The next result can be viewed as a
variation of their theorem. Although the result is not new, we
believe that it is a ``natural" way---from the mapping class groups
point of view---of thinking about the contact structures along with
their compatible open books. The reader may turn to \cite{k}, for
some facts we will use below about the mapping class groups.

{\Cor Any contact 3--manifold can be obtained by contact $(\pm
1)$--surgeries from a contact structure on some Seifert fibered
3--manifold.}

\begin{proof} Let $(Y, \xi)$ be a contact 3--manifold.
By Giroux's Theorem~\ref{giroux}, there is an open book $\OB$ of $Y$
compatible with $\xi$ whose page is a compact oriented surface
$\Sigma$ with $r$ boundary components, and whose monodromy is $\phi
: \Sigma \to \Sigma$. Let $\delta_i$ denote a curve parallel to the
$i$-th boundary component of $\Sigma$. Then there is a factorization
of $\phi$ into Dehn twists along some non-separating curves on
$\Sigma$ and curves parallel to the boundary components of $\Sigma$.
We can assume that $\phi= \phi^\prime \circ \psi $, where $\psi$
denotes the product of all the boundary parallel Dehn twists in the
factorization of $\psi$, since these curves are clearly disjoint
from the rest of the curves in the factorization. Moreover we can
also assume that $\psi = \prod_{i=1}^r t_{\delta_i}^{m_i}$ for some
nonzero integers $ m_1, m_2, \cdots, m_r$ by inserting some
canceling Dehn twists into the factorization of $\phi$ and using the
lantern relation in the mapping class groups. Here, if necessary, we
can increase the genus of the page by some positive stabilizations.
By Corollary~\ref{boundary}, there is a Seifert fibered 3--manifold
$Y_\psi$ with an open book $\OB_\psi$ whose page is $\Sigma$ and
whose monodromy is $\psi$. Now consider the contact structure
$\xi_\psi$ on $Y_\psi$ compatible with $\OB_\psi$. Embed all the
curves which appear in the factorization of $\phi^\prime $ on
different pages  of $\OB_\psi$. Since these curves are homologically
nontrivial on $\Sigma$, we can Legendrian realize these curves on
the convex pages of $\OB_\psi$ with respect to $\xi_\psi$. Now
applying contact $(-1)$--surgeries on curves which correspond to
right-handed Dehn twists and contact $(+1)$--surgeries which
correspond to left-handed Dehn twists yields $Y$ with the open book
$\OB$ and its compatible contact structure $\xi$. So we showed that
$(Y, \xi)$ can be obtained by contact $(\pm 1)$--surgeries along
Legendrian curves in the contact Seifert fibered 3--manifold
$(Y_\psi, \xi_\psi)$.

\end{proof}

\section{Surgery diagrams of some horizontal contact structures} \label{surg}

Consider a Seifert fibered 3--manifold with invariants
$$(g,n; -\frac{1}{p_1} , \cdots, -\frac{1}{p_k})$$ such that $n \leq
0 < p_i$, for $i=1,\cdots, k$. We can describe this 3--manifold as a
plumbing of circle bundles according to a star shaped graph with a
central vertex connected to $k$ other vertices. The central vertex
represents a bundle over a closed genus $g$ surface whose Euler
number is $n$ and the other $k$ vertices represent bundles over
$S^2$ with Euler numbers $p_i$, for $i=1,\cdots, k$.

We now show how to put this plumbing diagram into a standard form:
First observe that we can assume $p_i > 1$ for all $i$, since
otherwise we can blow down all the $(+1)$-vertices which will only
decrease $n$. Next blow up all the intersections of the central
vertex with the other $k$ vertices. The resulting plumbing graph
will have a central vertex with Euler number $n-k$ and $k$ linear
branches coming out of the central vertex whose vertices have Euler
numbers $(n-k, -1, p_i-1)$, including the central vertex, for
$i=1,\cdots, k$. Then blow up the intersection of $-1$ and $p_i-1$
at every branch to end up with the vertices having Euler numbers
$(n-k,-2, -1, p_i-2)$ on the $i$-th branch. Continue this process
$p_i-1$ times on the $i$-th branch, so that Euler numbers on each
branch becomes  $(n-k, -2,-2,\cdots, -2, -1, 1)$. Finally blow down
the last $+1$ to get the sequence $(n-k, -2,-2,\cdots, -2, -2)$, on
each branch, where the number of $-2$'s is equal to $p_i-1$ on the
$i$-th branch.

\begin{figure}[ht]

  \begin{center}

     \includegraphics{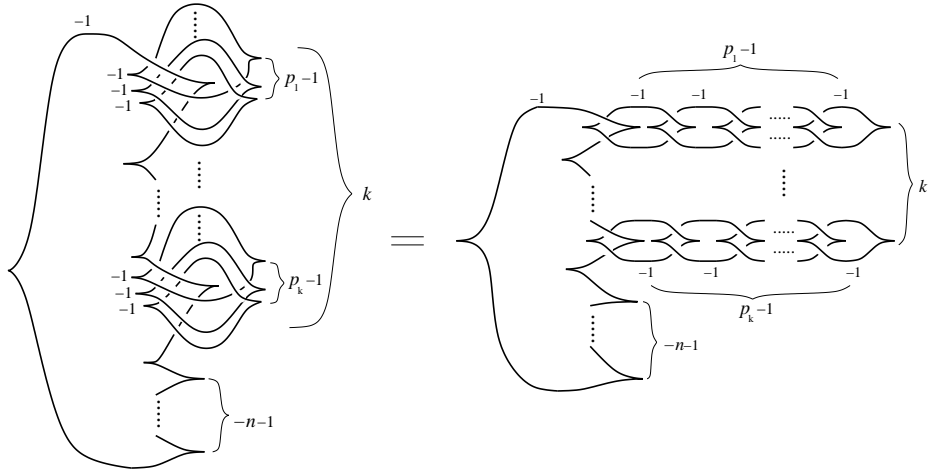}

   \caption{Horizontal contact structures on some Seifert fibered 3--manifolds}

\label{seifcont}

    \end{center}

  \end{figure}

Note that this is a non-positive standard plumbing
diagram for the given Seifert fibered 3--manifold. Now we observe
that the horizontal open book obtained by applying the algorithm in
\cite{eo} to this non-positive plumbing is isomorphic to the
horizontal open book we constructed in Section~\ref{construction}.
Suppose that $g=0$. Then by comparing this construction with
Sch\"{o}nenberger's method (cf. \cite{sc}) we conclude that the
horizontal contact structure compatible with the horizontal
open book on the Seifert fibered 3--manifold with invariants $(0,n;
-1/ p_1 , \cdots, -1/ p_k)$ where $n < 0 < p_i$ can be described by
the surgery diagram on the left-hand side in Figure~\ref{seifcont}.
Moreover one can show that the contact structure on the left-hand side
is isotopic to the contact structure depicted on the right-hand side.


\begin{thebibliography}{99999}

\bibitem[ACH]{ach}
C. Abbas, K. Cieliebak, H. Hofer, \emph{The Weinstein conjecture for
planar contact structures in dimension three};
arXiv:math.SG/0409355.

\bibitem[Al]{al}
J. Alexander, \emph{A lemma on systems of knotted curves,} Proc. Nat. Acad. Sci. USA
\textbf{9} (1923), 93--95.

\bibitem[DG]{dg}
F. Ding and H. Geiges, \emph{A Legendrian surgery presentation of
contact $3$-man\-i\-folds,} Math. Proc. Cambridge Philos. Soc.
{\bf136} (2004), 583--598.

\bibitem[EO]{eo}
T. Etg\"{u} and B. Ozbagci, {\em Explicit horizontal open books on
some plumbings}, to appear in the Inter. Jour. of Math.;
arXiv:math.GT/0509611


\bibitem[Et]{et}
J. Etnyre, {\em Lectures on open book decompositions and contact
structures}, Lecture notes from the Clay Mathematics Institute
Summer School on Floer Homology, Gauge Theory, and Low Dimensional
Topology at the Alfr\'{e}d R\'{e}nyi Institute;
arXiv:math.SG/0409402.

\bibitem[Gi]{gi}
E. Giroux, \emph{G\'{e}ometrie de contact: de la dimension trois
vers les dimensions sup\'{e}rieures,} Proceedings of the
International Congress of Mathematicians (Beijing 2002), Vol. II,
405--414.

\bibitem[HKM]{hkm}
K. Honda, W. Kazez, G. Matic, \emph{Right-veering diffeomorphisms of  compact surfaces with boundary I,} arXiv:math.GT/0510639

\bibitem[Ko]{k}
M. Korkmaz, \emph{Low-dimensional homology groups of mapping class
groups: a survey,} Turkish J.~Math. \textbf{26} (2002), 101--114.


\bibitem[LM]{lm}

P. Lisca and G. Matic, \emph{Transverse contact structures on Seifert 3--manifolds} Algebr. Geom. Topol. 4 (2004) 1125--1144.

\bibitem[Ma]{marti}
J. Martinet, \emph{Formes de contact sur les vari\'t\'es de dimension 3,} in: Proc.
Liverpool Singularity Sympos. II, Lecture Notes in Math. \textbf{209,} Springer, Berlin
(1971), 142--163.


\bibitem[OS]{ozst}
B. Ozbagci and A. Stipsicz, {\em Surgery on contact 3--manifolds and
Stein surfaces}, Bolyai Soc. Math. Stud., Vol. {\bf 13}, Springer,
2004.

\bibitem[Sc]{sc}
S. Sch\"{o}nenberger, {\em Planar open books and symplectic
fillings}, Ph.D. Dissertation, University of Pennsylvania, 2005.

\bibitem[TW]{tw}
W. Thurston and H. Winkelnkemper, \emph{On the existence of contact forms,} Proc. Amer.
Math. Soc. \textbf{52} (1975), 345--347.



\end{thebibliography}
\end{document}